\documentclass[11pt]{amsart}
\usepackage{url}
\usepackage{graphicx}
\usepackage{latexsym}
\usepackage{amsfonts,amsmath,amssymb}
\newtheorem{theorem}{Theorem}[section]																
\newtheorem{lemma}[theorem]{Lemma}

\usepackage{color}

\def\eps{\varepsilon}
\def\la{\lambda}
\def\a{\alpha}

\def\ga{\gamma}
\def\part{\partial}

\def\b1{\bold 1}

\newcommand{\beq}{\begin{equation}}
\newcommand{\eeq}{\end{equation}}

\theoremstyle{remark}

\numberwithin{equation}{section}
\date{\today}

\begin{document}

\title[Longest chains]{On increasing sequences formed by points from a random finite subset of a hypercube}
\author {Boris Pittel}
\address{Department of Mathematics, Ohio State, University, Columbus, OH 43210}
\email{pittel.1@osu.edu}

\maketitle

\begin{abstract} Consider $S$, a set of $n$ points chosen uniformly at random and independently from the unit hypercube
of dimension $t>2$. Order $S$ by using the Cartesian product of the $t$ standard orders of $[0,1]$. We determine a constant $\bar x(t)<e$ such that, with probability $\ge 1-\exp(-\Theta(\eps)n^{1/t})$, cardinality of a longest chain, i. e. a largest subset of comparable points, is at most $(\bar x(t)+\eps)n^{1/t}$. The bound $\bar x(t)$ complements an explicit lower bound obtained by Bollob\'as  and Winkler in 1988. Furthermore, we use Dilworth's theorem on partitions of a set into chains to prove that the cardinality of a largest antichain, i. e. a largest subset of incomparable points, is at least $(1-\eps) (n/e)^{1-1/t}$ with probability exponentially close to $1$.

\end{abstract}

\section{Introduction} Given two points $\bold X'$ and $\bold X''$ in $R^t$, we write $\bold X''\ge \bold X'$ if each component of $\bold X''$ is not strictly below the corresponding component of $\bold X'$. This makes $R^t$ a partially ordered set (poset).  Consider a set $\mathcal S_n=\{\bold X_i\}_{i\in [t]}$, where $\bold X_i$ are chosen independently, and uniformly, in the $t$-dimensional unit cube $\{\bold x=(x_1,\dots, x_t): x_i\in [0,1]\}$. Equivalently, the $nt$ components of $\bold X_1,\dots, \bold X_n$ are chosen uniformly, and independently of each other, in the interval 
$[0,1]$. The partial order on $R^t$ makes $\mathcal S_n$ a random poset.  

In this paper we study the length of the longest chain in $\mathcal S_n$, where a chain is defined as an ordered sequence 
$S=\{\bold X_{i_1}, \bold X_{i_2}, \cdots,\bold X_{i_{\ell}}\}$ such that $\bold X_{i_j}\le\bold X_{i_{j+1}}$, $1\le j\le \ell-1$; $\bold X_{i_0}:=\overrightarrow{\bold 0}=(0,\dots, 0)$.
In 1961, Ulam \cite{Ula} posed a problem to determine an asymptotic length of the longest subsequence in the uniformly random
permutation of $[n]$. It is easy to see that Ulam's problem is equivalent to the problem of the longest chain in the case $t=2$. 

In 1977, the results of Logan and Shepp \cite{Log} and Vershik and Kerov \cite{Ver1}, \cite{Ver2} provided a solution: with high probability, Ulam's parameter is asymptotic to $2n^{1/2}$. This breakthrough result was based on asymptotic analysis of 
the shape of a Young diagram of size $n$, chosen according to Plancherel distribution: the probability of a generic diagram $\la$ is $d^2(\la)/n!$, where $d(\la)$ is $n!$ divided by the product of $\la$'s hook lengths, which happens to be the number of the Young tableaux of shape $\la$. Combinatorially,
the distribution arises from a Robinson--Schensted bijection (see \cite{Sch}, Knuth \cite{Knu}, Vol. 3, Theorem A), between $n!$ permutations and pairs of Young tableaux of the same shape $\la$, with a property that the size of the longest increasing subsequence in a permutation equals the size of the longest row in the associated $\la$. In 1991, Aldous and Diaconis \cite{Ald} found an alternative proof of the $2n^{1/2}$ result based on Hammersley's interacting particles process.

Prior to \cite{Log}, \cite{Ver1}, and \cite{Ver2}, Hammersley \cite{Ham1}, \cite{Ham2} was already able to prove that the parameter in question scaled by $n^{1/2}$
converges, in probability, to a constant sandwiched between $\pi/8$ and $e$. That $e$ is an upper bound for the limit 
in question followed from an observation that the expected number of chains of length $\ell$ is $(n)_{\ell}/(\ell!)^2$. The lower bound followed from an asymptotic analysis of an auxiliary process of growing a chain of points via a two-dimensional Poisson process. Convergence of the scaled length of a longest subsequence to a constant limit followed from a general result for subadditive ergodic processes, due to Kingman \cite{Kin}. 

Sixteen years later, Bollob\'as and Winkler \cite{Bol} studied the case $t>2$. 
Extending Hammersley-Kingman's techniques, they proved that the length of the longest chain scaled by $n^{1/t}$ converges in probability to a constant in  
$\bigl[t^2/t!^{1/t}\Gamma(1/t),e\bigr)$, the left endpoint converging to the (excluded) right endpoint $e$ as $t\to\infty$. 
In this paper we replace the upper bound $e$ with a smaller, $t$-dependent, constant. 

Associated with a chain $S$ there are $\ell$ {\it intervals\/}: 
\[
[\overrightarrow{\bold 0},\bold X_{i_1}],\, [\bold X_{i_1},\bold X_{i_2}], \dots, [\bold X_{i_{\ell-1}},
\bold X_{i_{\ell}}],
\]
each interval consisting of the points $\bold X_j$, ($j\neq i_1,\dots, i_{\ell}$), comparable to the interval's endpoints, and sandwiched strictly between them. We call a chain $\{\bold X_{i_1}\le \bold X_{i_2}\le \cdots \le \bold X_{i_{\ell}}\}$
{\it maximal\/} if each of the $\ell$ intervals is empty. Clearly, the longest chains are maximal. (Of course, for a longest chain $\bold X_{i_1}\le\cdots\le \bold X_{i_{\ell}}$, the extra interval $[\bold X_{i_{\ell}},\overrightarrow{\bold 1}]$,
$\overrightarrow{\bold 1}:=(1,\dots, 1)$, is empty too.) 
\begin{theorem}\label{lem1} Let $C_{\ell}$ denote the total number of maximal chains of length $\ell=x n^{1/t}$, $x\le e$. {\bf (i)\/} For $a>\ga t$, ($\ga=-\Gamma'(1)$), and $b>0$, introduce 
\[
q(x;a,b)\!=\!\max\!\left\{\!\min_{u\in (0,1)}\frac{\Gamma^t(1-u)}{e^{au}}; \frac{1+b}{e^b};\exp\left(-\frac{e^{-a}}
{x^t(1+b)^t}\right)\!\right\}\!<1.
\]
 $q(x)=\min_{a,b}q(x;a,b)$ is attained at a unique point $(a(x), b(x))$ such that 
 \[
q(x)=\min_{u\in (0,1)}\frac{\Gamma^t(1-u)}{e^{a(x)u}}= \frac{1+b(x)}{e^{b(x)}}=\exp\left(-\frac{e^{-a(x)}}
{x^t(1+b(x))^t}\right)<1;
\]
$q(x)$ decreases with $x$.  {\bf (ii)\/} Define $\bar x=\inf\bigl\{x\in [e^{-\ga},e]: \bigl(\frac{e}{x}\bigr)^tq(x)<1\bigr\}$. $\bar x<e$ since $q(e)<1$, and
 \[
 \Bbb E[C_{\ell}]=O\bigl[\exp(-\Theta(\eps)n^{1/t})\bigr], \text{ if }  \ell =(\bar x+\eps)n^{1/t}, \eps>0. 
\]
 Consequently, the length of the longest chain is below $(\bar x+\eps) n^{1/t}$ with probability $\ge 1-O\bigl[\exp(-\Theta(\eps)n^{1/t})\bigr]$.
\end{theorem}
\noindent {\it Note.\/} Of interest in its own right, the characterization of $q(x)$ may well be instrumental for a numerical analysis of $\bar x(t)$.

Our second result is about the size of the largest antichain in $\mathcal S_n$, i.e. the size of the largest subset of $\mathcal S_n$
consisting of pairwise incomparable points $\bold X_i$.

\begin{theorem}\label{thm2} Let $D_n$ be size of the largest antichain in $\mathcal S_n$. Then
\[
\Bbb P\left(D_n\ge (1-\eps) (en)^{1-1/t}\right)\ge 1- \rho^n(\eps),\quad \rho(\eps) \in (0,1).
\]
\end{theorem}
\noindent The proof is based on Dilworth's theorem, stating that, for every partially ordered set, the maximum size of an antichain equals the minimum number of chains in a partition of the set into chains. 

In conclusion, we refer the reader to Steele \cite{Ste} for a highly readable, meticulously written survey of a rich collection of papers on lengths of increasing subsequences in a variety of combinatorial settings.

\section{Proof of Theorem \ref{lem1}} 
 
By symmetry, with $\bold x_0=\overrightarrow{\bold 0}$,
 \[
 \Bbb E[C_{\ell}]=(n)_{\ell}\int\limits_{\bold 0\le \bold x(1)\le \bold x_2\le\dots\le \bold x_{\ell}\le \bold 1}\!\!\left(\!1-\sum_{i=1}^{\ell}
 \prod_{j=1}^t\bigl(x_j(i)-x_j(i-1)\bigr)\!\right)^{n-\ell}\!\!\prod_{i=1}^l d\bold x(i).
 \]
{\it Explanation\/}:   $(n)_{\ell}$ is the total number of ways to select an ordered sequence of points $\{\bold X(i)\}_{i\in [\ell]}$ in the cube $[0,1]^t$. The integral is the probability that this sequence is increasing {\it and\/} that all the $\ell$ intervals $[\bold X(i-1),\bold X(i+1)]$ are empty; $\bold X(0):=\overrightarrow{\bold 0}$.  Indeed, conditioned on 
$\{\bold X(i)=\bold x(i)\}_{i\in [\ell]}$, the integrand is the probability that none of the remaining $(n-\ell)$ points $\bold X(k)$ belongs to any of the $\ell$ intervals in question.

Introduce the $\ell$ increments $\bold y(i)=\bold x(i)-\bold x(i-1)$, $i\in [1,\ell]$, so that $\sum_{i=1}^{\ell} \bold y(i)\le \bold 1$.
Switching from the variables $\bold x(i)$, $i\in [\ell]$, to $\bold y(i)=\{y_j(i)\}_{j\in [t]}$, $i\in [\ell]$, we rewrite the above formula:
\[
\Bbb E[C_{\ell}]=(n)_{\ell}\int\limits_{\bold y_1,\dots,\bold y_{\ell}\ge \bold 0,\atop \bold y_1+\cdots+\bold y_\ell\le \bold 1}
\biggl(1-\sum_{i\in [\ell]}\prod_{j\in [t]}y_j(i)\biggr)^{n-\ell}\,\prod_{i=1}^{\ell} d\bold y(i).
\]
Our task is to upper bound this $t\ell$-dimensional integral. As the first step, we use the arithmetic/geometric means inequality in combination with $1-z\le e^{-z}$, $(z\in [0,1])$, and bound the integrand by
\[
\biggl(1-\ell\prod_{i\in [\ell]}\prod_{j\in [t]}y_j^{1/\ell}(i)\biggr)^{n-\ell}
\le\exp\biggl(\!\!-\ell(n-\ell)\prod_{j\in [t]}\prod_{i\in [\ell]}y_j^{1/\ell}(i)\biggr).
\]
Probabilistically, 
\begin{equation}\label{12.9}
\Bbb E[C_{\ell}]\le (n)_{\ell}\,\Bbb E\biggl[\exp\biggl(-\ell(n-\ell)\prod_{j\in [t]}\prod_{i\in [\ell]}Y_j^{1/\ell}(i)\biggr)
\cdot\Bbb I\biggl(\sum_{i\in [\ell]}\bold Y(i)\le \overrightarrow{\bold 1}\biggr)\biggr],
\end{equation}
where $t\ell$ random variables $\{Y_j(i)\}_{i\in [\ell], j\in [t]}$ are independent $(0,1)$-Uniforms, and $\bold Y(i)=\{Y_j(i)\}_{j\in [t]}$.

To proceed, for each $j\in [t]$ switch to $s_j=\sum_{i\in [\ell]}y_j(i)$, $v_j(i)=y_j(i)/s_j$, $i<\ell$, and set $v_j(\ell)=1-\sum_{i<\ell}v_j(i)$. 
In the new variables, 
\[
\prod_{j\in [t]}\prod_{i\in [\ell]}y_j^{1/\ell}(i)=\Pi(\bold v) \cdot\prod_{j\in [t]} s_j, \quad \Pi(\bold v)=\prod_{j\in [t]}\prod_{i\in [\ell]}(v_j(i))
^{1/\ell}.
\]
An admissible $\{s_j, v_j(i)\}_{j\in [t], i\in [\ell]}$ satisfies the conditions $0\le v_j(i)\le \min \{1,1/s_j\}$, $s_j\le 1$. The Jacobian of 
$\{y_j(i)\}_{i\in [\ell], j\in [t]}$ with respect to $\{s_j, v_j(i)\}_{i\in [\ell-1],\,j\in [t]}$ is $\prod_{j\in [t]} s_j^{\ell-1}$, implying  that the joint density of $\{S_j:=\sum_{i\in [\ell]}Y_j(i), V_j(i):=Y_j(i)/S_j\}_{j\in [t], i\in [\ell-1]}$ is {\it at most\/}
\[
\prod_{j\in [t]}s_j^{\ell-1}\cdot\Bbb I\biggl(\sum_{i\in [\ell-1]}v_j(i)\le 1\biggr)=\prod_{j\in [t]}\frac{s_j^{\ell-1}}{(\ell-1)!}\cdot g\bigl(\{v_j(i)\}_{i\in [\ell-1]}\bigr).
\]
Here, for each $j\in [t]$, 
\[
g(\{v_j(i)\}_{i\in [\ell-1]}):=(\ell-1)!\cdot \Bbb I\biggl(\sum_{i\in [\ell-1]}v_j(i)\le 1\biggr)
\]
is the joint density of the {\it first\/} $(\ell-1)$ components of $\bold L_j=\{L_j(i)\}_{i\in [\ell]}$, where $L_j(i)$ are the lengths of 
$\ell$ consecutive subintervals of $[0,1]$ obtained by selecting $(\ell-1)$ points independently and uniformly at random in $[0,1]$. The product of these densities signifies that $\ell$-dimensional vectors $\bold L_1,\dots, \bold L_t$ are independent.
So, the expectation on the the RHS of \eqref{12.9} is at most
\begin{multline}\label{13}
\frac{1}{((\ell-1)!)^t}\int\limits_{\{v_j(i)\}}\prod_{j\in [t]} g\bigl(\{v_j(i)\}_{i\in [\ell-1]}\bigr)\\
\times\left(\,\int\limits_{\{s_j\}} \exp\biggl(-\ell(n-\ell)\cdot\Pi(\bold v)\prod_{j\in [t]}s_j\biggr)\prod_{j\in [t]}\!s_j^{\ell-1}\,\,
d\bold s\right)\,\prod_{j\in [t]} \prod_{i\in [\ell-1]} d\bold v_j(i);
\end{multline}
here $\sum_{i\in[\ell]}v_j(i)\le 1$, $s_j\le 1$, ($j\in [t]$). 

To evaluate the innermost integral, switch to the variables $z_1=s_1, z_2=s_1s_2,\dots, z_{t-1}=s_1\cdots s_{t-1}, z=
s_1\cdots s_t$. The new variables meet the conditions $1\ge z_1\ge\cdots\ge z_{t-1}\ge z$. The Jacobian of the variables $s_j$, with respect to $z_1,\dots, z_{t-1}, z$ is $1/z_1\times\dots\times 1/z_{t-1}$. So, the new integrand is 
\[
z^{\ell-1}\exp\biggl(-\ell(n-\ell)\cdot\Pi(\bold v)z\biggr) \cdot\prod_{j\in [t-1]}z_j^{-1}.
\]
Integrating this function with respect to $z_1,\dots, z_{t-1}$, first over $z_1\in [z_2,1]$, second over 
$z_2\in [z_3,1]$,..., and last over $z_{t-1}\in [z,1]$, we see that the innermost integral simplifies to
\[
\frac{1}{(t-1)!}\int_0^1z^{\ell-1}\exp\left[-\ell(n-\ell)\cdot \Pi(\bold v) z\right] \log^{t-1}(1/z) \,dz.
\]
Hence, recalling the definition of $g(\{v(i)\}_{i\in [\ell-1]})$, and setting $\bold L:=\{\bold L_j\}_{j\in [t]}$, we rewrite \eqref{13} this
way:
\begin{equation*}
\frac{\ell^t}{((\ell)!)^t\,(t-1)!}\,\Bbb E\left[\int_0^1z^{\ell-1}\exp\left[-\ell(n-\ell)\cdot \Pi(\bold L) z\right] \log^{t-1}(1/z \,dz\right].
\end{equation*}
In summary, \eqref{12.9} becomes
\begin{equation}\label{14}
\begin{aligned}
\Bbb E[C_{\ell}]&\le \frac{(n)_{\ell}\,\ell^t}{(\ell !)^t (t-1)!}\\
&\quad\times\Bbb E\left[\int_0^1z^{\ell-1}\exp\left[-\ell(n-\ell)\cdot \Pi(\bold L) z\right] \log^{t-1}(1/z \,dz\right].
\end{aligned}
\end{equation}
To sharply estimate the RHS expectation, we use a classic property of the uniformly random partition of $[0,1]$,
(Karlin and Taylor \cite{Kar}): 
the lengths $L(1),\dots, L(\ell)$ of the subintervals of $[0,1]$ obtained by throwing $\ell-1$ points into $[0,1]$ uniformly and independently have the same joint distribution as 
$\bigl\{W(1)/S(\ell),\dots, W(\ell)/S(\ell)\bigr\}$, where $W(i)$ are independent exponentials with mean $1$, and $S(\ell)=\sum_{i\in [\ell]}W(i)$. Then, introducing $t\ell$ independent exponentials $W_j(i)$, and using
$\overset{\mathcal D}\equiv$ to indicate equality of two distributions in question, we have
\[
\Pi(\bold L)\overset{\mathcal D}\equiv\Pi(\bold W)=\frac{\prod_{j,i}W^{1/\ell}_j(i)}{\prod_jS_j(\ell)},\quad S_j(\ell)=\sum_iW_j(i).
\]
{\bf (1)\/} Using Chernoff-type bound for the i.i.d. summands $\log W_j(i)$, we have: for $a>0$,
\begin{multline*}
\Bbb P\biggl(\ell^{-1}\!\!\!\!\sum_{j\in [t], i\in [\ell]}\!\!\! \log W_j(i)\le -a\biggr)=
\Bbb P\biggl(\sum_{j\in [t], i\in [\ell]}\!\!\! \log W_j(i)\le -a\ell\biggr)\\
\le\inf_{u\in (0,1)}\frac{\Bbb E\biggl[\exp\biggl(-u\sum_{j\in [t], i\in [\ell]}\log W_j(i)\biggr)\biggr]}{\exp(a\ell u)}
= \inf_{u\in (0,1)}\frac{\bigl(\Bbb E[e^{-u\log W}]\bigr)^{t\ell}}{e^{a\ell u}};
\end{multline*}
here 
\[
\Bbb E[e^{-u\log W}]=\int_0^{\infty}e^{-u\log w} e^{-w}\,dw=\int_0^{\infty} w^{-u} e^{-w}\,dw=\Gamma(1-u).
\]
Therefore
\begin{equation*}
\Bbb P\biggl(\ell^{-1}\!\!\!\!\sum_{j\in [t], i\in [\ell]}\!\!\! \log W_j(i)\le -a\biggr)\le \left(\inf_{u\in (0,1)}\frac{\Gamma^t(1-u)}{e^{au}}\right)^{\ell}.
\end{equation*} 
Since $\Gamma'(1)=-\ga$, where $\ga$ (Euler-Mascheroni constant) is $0.57721\dots$), the infimum is strictly below $1$ if $a>\ga t$, which we assume from now. Since $\Gamma(x)$ is log-convex (Andrews, Askey, and Roy \cite{And}), and $\Gamma(+0)=\infty$, the infimum is attained at a single stationary point $u(a)$ of $e^{-au}\, \Gamma^t(1-u)$ in $(0,1)$, i.e. the $a$-dependent root of 
\[
\bigl(t\log\Gamma(1-u)-au\bigr)'=-\left.t\frac{\Gamma'(v)}{\Gamma(v)}\right|_{v=1-u}-a=0.
\]
Hence
\begin{equation}\label{15.5}
\Bbb P\biggl(\ell^{-1}\!\!\!\!\sum_{j\in [t], i\in [\ell]}\!\!\! \log W_j(i)\le -a\biggr)\le \rho^{\ell}(a),\quad \rho(a)
:=\min_{u\in (0,1)}\frac{\Gamma^t(1-u)}{e^{au}}.
\end{equation}
$\rho(a)$ decreases as a function of $a$, $\rho(\ga t+)=1$, and $\rho(\infty)=0$.


{\bf (2)\/} Furthermore, $S_j(\ell)$, $j\in [t]$, are i.i.d., with density $e^{-y} y^{\ell-1}/\Gamma(\ell)!$. Picking $b>0$, we bound, again via Chernoff-type method, 
\begin{multline*}
\Bbb P\biggl(\sum_{j\in [t]}\log S_j(\ell)\ge t\log(\ell(1+b))\biggr)\le \frac{\Bbb E^t\bigl[e^{u\log S_1(\ell)}\bigr]}{\exp\bigl[ut\log(\ell(1+b)\bigr])}\\
=(\ell(1+b))^{-ut}\left(\int_0^{\infty}e^{-y}\frac{y^{\ell-1}}{\Gamma(\ell)} y^u\,dy\right)^t\\
=\left[(\ell(1+b))^{-u}\frac{\Gamma(u+\ell)}{\Gamma(\ell)}\right]^t.
\end{multline*}
By Stirling-type formula $\Gamma(\eta)=\Theta\bigl(\eta^{1/2}(\eta/e)^{\eta}\bigr)$, $u=\ell b$ delivers qualitatively best bound,
namely 
\begin{equation}\label{16}
\Bbb P\biggl(\sum_j\log S_j(\ell)\ge t\log(\ell(1+b))\biggr)
\le \left(\frac{1+b}{e^b}\right)^{\ell}.
\end{equation}

Combining \eqref{15.5} and \eqref{16}, we bound the expectation in \eqref{14} by 
\begin{multline}\label{17}
O\biggl[\rho^{\ell}(t)+\left(\frac{1+b}{e^b}\right)^{\ell}\biggr] \int_0^1 z^{\ell-1}\log^{t-1}(1/z)\,dz\\
+\int_0^1z^{\ell-1}\exp\left(\!-e^{-a}\frac{\ell(n-\ell)}{\ell^t(1+b)^t}\,z\!\right)\cdot\log^{t-1}(1/z)\,dz.
\end{multline}
Integrating the top integral in \eqref{17} by parts once, we obtain a recurrence equation, which implies that the integral equals $(t-1)!/\ell^t$. Since $t\ge 2$, the bottom integral is of order
\begin{multline*}
\int_0^1z^{\ell-1}\exp\left(\!-e^{-a}\frac{\ell n}{\ell^t(1+b)^t}\,z\!\right)\cdot\log^{t-1}(1/z)\,dz\\
\le \max_{z\in [0,1]} z^{\ell-1}\exp\left(\!-e^{-a}\frac{\ell n}{\ell^t(1+b)^t}\,z\!\right)
\times\int_0^1\log^{t-1}(1/z)\,dz.
\end{multline*}
The second line integral equals $(t-1)!$, and the outside maximum is attained at $z=1$, provided that $\ell\ge e^{-\ga} n^{1/t}$ and $n$ is large. Indeed, the function  in question increases for $z\le (1-1/\ell)e^{a}(1+b)^t (\ell^t/n)$; since $a> \ga t$, the bound for $z$ exceeds $1$ for $\ell\ge e^{-\ga} n^{1/t}$, and $n$ is sufficiently large.
The resulting bound together with \eqref{17} imply that, for $x:=\ell/n^{1/t}\ge e^{-\ga}$, the expectation in \eqref{14} is of order
\begin{multline*}
\rho^{\ell}(a)+\left(\frac{1+b}{e^b}\right)^{\ell}+\exp^{\ell}\left(-\frac{e^{-a}}{x^t(1+b)^t}\right)
\le\rho^{\ell}(a)+2\left(\frac{1+b(a,x)}{e^{b(a,x)}}\right)^{\ell},
\end{multline*}
where $b(a,x)$ is a unique root of
\begin{equation}\label{18}
\frac{1+b}{e^b}=\exp\left(-\frac{e^{-a}}{x^t(1+b)^t}\right),\quad b>0.
\end{equation}  
And $b(a,x)$ is the minimum point of $\max\left\{\frac{1+b}{e^{b}}, \exp\left(-\frac{e^{-a}}{x^t(1+b)^t}\right)\right\}$,
since as functions of $b$, $\frac{1+b}{e^b}$ decreases from $1$ to $0$, while $\exp\left(-\frac{e^{-a}}{x^t(1+b)^t}\right)$ increases to $1$.
Now, for 
\[
H(x,a,b):=\log(1+b)-b+\frac{e^{-a}}{x^t(1+b)^t}, 
\]
we have $H'_b<0$, $H'_a<0$, and $H'_x<0$. 
By implicit differentiation, $b(a,x)$ decreases as a function of $a$ and $x$.
Finally, we choose $a=a(x)>\ga t$ such that 
\begin{equation}\label{19}
\rho(a)=\frac{1+b(a,x)}{e^{b(a,x)}};
\end{equation}
$a(x)$ exists since, as functions of $a\in (\ga t,\infty)$, $\rho(a)$ decreases from $1$ to $0$, while $\frac{1+b(a,x)}{e^{b(a,x)}}$ increases to $1$. And $a(x)$ minimizes $\max\left\{\rho(a), \exp\left(-\frac{e^{-a}}{x^t(1+b)^t}\right)\right\}$.
Now, for
\[
K(x,a, b(x,a)):=\rho(a)-\frac{1+b(a,x)}{e^{b(a,x)}},
\]
we have 
\begin{multline*}
K'_a(x,a,b(x,a))=K'_a(x,a, u)\big|_{u=b(x,a)}+K'_u(x,a,u)\big|_{u=b(a,x)}b'_a(a,x)\\
=\rho'_a(a)+b(x,a)e^{-b(x,a)}b'_a(a,x)<0,
\end{multline*}
since $\rho'_a(a)<0$ and $b'_a(a,x)<0$, and 
\begin{multline*}
K'_x(x,a,b(x,a))=K'_x(x,a,u)\big|_{u=b(x,a)}+K'_u(x,a,u)\big|_{u=b(a,x)}b'_x(a,x)\\
=b(x,a)e^{-b(x,a)}b'_x(a,x)<0.
\end{multline*}
Again by implicit differentiation, $a(x)$, the root of $K(x,a, b(x,a))=0$,  decreases as a function of $x$

In summary, we proved that $(a(x), b(x))$ is a unique minimum point of 
\begin{equation*}
q(x;a,b):=\max\left\{\rho(a), \frac{1+b}{e^b}, \exp\biggl(-\frac{e^{-a}}{x^t(1+b)^t}\biggr)\right\}, 
\end{equation*}
and, denoting $q(x)=q(x; a(x), b(x))$, 
\[
q(x)=\rho(a(x))=\frac{1+b(x)}{e^{b(x)}}=\exp\biggl(-\frac{e^{-a(x)}}{x^t(1+b(x))^t}\biggr).
\]
Since both $\rho(a)$ and $a(x)$ are decreasing, $q(x)=\rho(a(x))$ is increasing.

Now, in \eqref{14},  $(n)_{\ell} \ell^t/(\ell!)^t$ is of order $\bigl[n/(\ell/e)^t\bigr]^{\ell}\ell^{t/2}$. So, the previous discussion shows that for $\ell=x n^{1/t}$, $E[C_{\ell}]$ is at most of order $n^{1/2}\left(\left(\frac{e}{x}\right)^tq(x)\right)^{\ell}$.
Introducing $\bar x=\inf\{x\in [e^{-\ga},e]: (e/x)^tq(x)<1\}$,  we see that $\bar x<e$, since $q(e)<1$.
And we conclude that $\Bbb E[C_{\ell}]=O\bigl[\exp(-\Theta(\eps)n^{1/t})\bigr]$, if $\ell\ge (\bar x+\eps)n^{1/t}$.
The proof of Theorem \ref{lem1} is complete.

\section{Proof of Theorem \ref{thm2}}

By a general theorem, due to Dilworth (Anderson \cite{Ande}) $D_n$, the size of the largest multichain in $\mathcal S_n$ is the 
smallest number of disjoint chains needed to cover all points of $\mathcal S_n$. 
\begin{lemma}\label{lem2}
Introducing $N(n,k)$, the total number of collections of $k$ point-disjoint chains that collectively cover $\mathcal S_n$, we have
\[
\Bbb E[N(n,k)]=\sum_{\{k_j\}_{j\ge 1}\atop \sum_j jk_j=n,\, \sum_jk_j=k}\frac{n!}{\prod_jk_j! (j!)^{k_j}}\cdot \prod_{j'}\left(\frac{j'!}{((j'!)^t}\right)^{k_{j'}}.
\]
\end{lemma}
\begin{proof} The first factor is the total number of ways to partition $[n]$ into $k$ subsets, such that  the number of subsets with $j$ elements equals $k_j$. The second factor is the probability that for each of those partitions the points $\bold X_i$, with $i$ belonging to every one of the subsets, are all comparable to each other, i.e. form a chain. 

Indeed, if such a subset has $j$ elements $i_1,\dots, i_j$, then all $(j!)^t$ ways to order the first components, the second components, and the last components of
$\bold X_{i_1},\dots, \bold X_{i_j}$ are equally likely, and only $j!$ ways for all $t$ of those components-confined orders to coincide.  
\end{proof}

By this Lemma,
\begin{multline*}
\sum_{n,k}\frac{x^n y^k}{n!}\Bbb E[N(n,k)]=\sum_{\{k_j\}_{j\ge 1}}\prod_{j\ge 1}\frac{(x^jy/(j!)^t)^{k_j}}{k_j!}\\
=\prod_{j\ge 1}\exp\left(\frac{yx^j}{(j!)^t}\right)=F(x,y):=\exp\left(y\sum_{j\ge 1}\frac{x^j}{(j!)^t}\right),
\end{multline*}
an explicit exponential generating function for the bivariate sequence $\{\Bbb E[N(n,k)]\}$. 

So, invoking Dilworth's theorem, and following it with a Chernoff-type bound, for each $d$, $x>0$ and $y\in (0,1)$, we have
\begin{multline*}
\Bbb P(D_n\le d)\le \sum_{k\le d}\Bbb E[N(n,k)]=n!\sum_{k\le d}\frac{\Bbb E[N(n,k)]}{n!}\\
\le n!\, y^{-d}\sum_{k\ge 1}\frac{\Bbb E[N(n,k)]}{n!}y^k\le n!\, x^{-n}y^{-d}\sum_{\nu,k}\frac{\Bbb E[N(\nu,k)]}{\nu!} x^{\nu} y^k\\
=n! \,x^{-n}y^{-d} F(x,y).
\end{multline*}
To upper bound $F(x,y)$ without sacrificing much of accuracy, we bound 
\begin{multline*}
\sum_{j\ge 1}\frac{x^j}{(j!)^t}=\sum_{j\ge 1}\frac{(jt)!}{(j!)^t}\cdot \frac{x^j}{(jt)!}\le \sum_{j\ge 0} t^{jt}\cdot\frac{(x^{1/t})^{jt}}{(jt)!}
\le \sum_{\nu\ge 0}\frac{(t x^{1/t})^{\nu}}{\nu!}=e^{tx^{1/t}}.\\
\end{multline*}
Therefore
\[
\Bbb P(D_n\le d) \le n!\,x^{-n}y^{-d}\exp\bigl(y e^{tx^{1/t}}\bigr), \quad \forall\, x>0,\, y\in (0,1).
\]
Further, the RHS function of $x,y$ has a unique stationary point $(x (n,d),y (n,d))$ given by
\[
x(n,d)=\left(\frac{n}{d}\right)^t,\quad y(n,d)=d e^{-tn/d}.
\]
The function $f(\eta)=\eta e^{-tn/\eta}$ is increasing, and $f(n/\log n)=O(n^{-(t-1)})\to 0$. Thus $y(n,d)<1$ for $d\le n/\log n$
and $n$ sufficiently large. So, for $d\le n/\log^2 n$,
\begin{multline*}
\Bbb P(D_n\le d) \le n!\,x^{-n}(n,d)y^{-d}(n,d)\exp\bigl(y(n,d) e^{ty(n,d)x^{1/t}(n,d)}\bigr)\\ 
=n!\left(\frac{n}{d}\right)^{-tn}\bigl(d e^{-tn/d}\bigr)^{-d}\cdot\exp\bigl(ne^t e^{-tn/d}\bigr)\\
=\left[\frac{n}{e}\cdot\left(\frac{n}{ed}\right)^{-t}\cdot (1+O(\log^{-1}n))\right]^n.\\
\end{multline*}
In particular, for $d=\a n^{1-1/t}$,
\[
\Bbb P(D_n\le \a n^{1-1/t})\le\left[e^{-1}(e\a)^t\cdot(1+O(\log^{-1}n))\right]^n,
\] 
and $e^{-1}(e\a)^t<1$ if $\a<e^{1/t-1}$. Therefore, if $\eps\in (0,1)$, then 
\[
\Bbb P\left(D_n\le (1-\eps)\left(\frac{n}{e}\right)^{1-1/t}\right)\le q^n(\eps),\quad q(\eps)\in (0,1).
\]
The proof of Theorem \ref{thm2} is complete.


\begin{thebibliography}{99}
\bibitem{Ald}
D. Aldous and P. Diaconis, \textit{Hammersley's interacting particle process and longest increasing subsequences}, Technical Report, Dept. Stat., U.C. Berkeley, CA (1993).
\bibitem{Ande}
I. Anderson, \textit{Combinatorics of Finite Sets}, Oxford (1987).
\bibitem{And}
G. E. Andrews, R. Askey, and R. Roy, \textit{Special functions}, Encyclopedia Math. Appl. Cambridge (1999).
\bibitem{Bol}
B. Bollob\'as and P. Winkler, \textit{Chain among random points in Euclidean space}, Proc, Amer. Math. Soc.
{\bf 103} (1988), 347--353.
\bibitem{Log}
B. F. Logan and L. A. Shepp, \textit{A variational problem for random Young tableaux}, Adv. in Math. {\bf 26} (1977), 
206--222.
\bibitem{Ham1}
J. M. Hammersley, \textit{A few seedlings of research}, Proc. Sixth Berkeley Sympos. Math. Stat. Prob., Univ. of California Press, Berkeley (1972), 345--394.
\bibitem{Ham2}
-----------------, \textit{Postulates for subadditive processes}, Ann. Probab. {\bf 2} (1973), 903.
\bibitem{Kar}
S. Karlin and H. M. Taylor, \textit{A second course in stochastic processes}, Acad. Press (1981).
\bibitem{Kin}
J. F. C. Kingman, \textit{Subadditive ergodic theory}, Ann. Probab. {\bf 1} (1973), 883--909.
\bibitem{Knu}
D. E. Knuth, \textit{The Art of Computer Programming}, 
Vol 3, Sorting and Searching, 2nd Edition, Section 5.2.2. 
\bibitem{Sch}
C. Schensted, \textit{Longest increasing and decreasing subsequences}, Canad, J. Math. {\bf 13} (1961), 179--191.
\bibitem{Ste}
J. M. Steele, \textit{Variations on the monotone subsequence theme of Erd\"os and Szekeres}, Discrete probability and algorithms, IMA Vol. Math. Appl., Springer-Verlag {\bf 72} (1993) 111--131.
\bibitem{Ula} 
S. M. Ulam, \textit{Monte Carlo calculations in problems of mathematical physics}, Modern Mathematics for the Engineer (E. F. Beckenbach, Ed.),  McGraw-Hill (1961).
\bibitem{Ver1}
A. M. Versik and S. V. Kerov, \textit{Asymptotics of the Plancherel measure of the symmetric group and the limiting form of Young tableaux}, Dokl. Akad. Nauk SSSR {\bf 233} (1977), 1024--1028.
\bibitem{Ver2}
-----------------, \textit{Asymptotics of the largest and the typical dimensions of irreducible representations of the symmetric group}, Funct. Anal. Appl. {\bf 19} (1985), 21--31.

\end{thebibliography}
\end{document}